\documentclass[leqno]{amsart}
\usepackage{amsmath,amsthm,amsfonts,amssymb}
\usepackage[backref=page]{hyperref}
\hypersetup{colorlinks=true,citecolor=blue,linkcolor=blue,urlcolor=blue,pdfstartview=FitH,
pdfauthor=Valentin Blomer ,pdftitle=Subconvexity for a double Dirichlet series}

\theoremstyle{plain}
\newtheorem{satz}{Theorem}
\newtheorem{lemma}{Lemma}

\theoremstyle{definition}

\renewcommand{\geq}{\geqslant}
\renewcommand{\leq}{\leqslant}

\begin{document}

\title{Subconvexity for a double Dirichlet series}

\author{Valentin Blomer}
\address{Mathematisches Institut, Bunsenstr. 3-5, 37073 G\"ottingen} \email{blomer@uni-math.gwdg.de} 
\thanks{author supported in part by the Volkswagen Foundation and a Sloan Research Fellowship}


\begin{abstract}  For two real characters $\psi, \psi'$ of conductor dividing 8 define
\begin{displaymath}
  Z(s, w; \psi, \psi') := \zeta_2(2s+2w-1) \sum_{d \text{ odd}} \frac{L_2(s, \chi_d\psi) \psi'(d)}{d^w}
\end{displaymath}
where $\chi_d = \left(\frac{d}{.}\right)$ and the subscript 2 denotes that the Euler factor at 2 has been removed. These double Dirichlet series can be extended to $\Bbb{C}^2$ possessing a group of functional equations isomorphic to  $D_{12}$. The convexity bound for $Z(s, w; \psi, \psi')$ is $|s w(s+w)|^{1/4+\varepsilon}$ for $\Re s = \Re w = 1/2$. It is proved that
\begin{displaymath}
  Z(s, w; \psi, \psi') \ll |s w(s+w)|^{1/6+\varepsilon}, \quad \Re s = \Re w = 1/2. 
\end{displaymath}
Moreover, the following mean square Lindel\"of type bound holds
\begin{displaymath}
 \int_{-Y_1}^{Y_1} \int_{-Y_2}^{Y_2} |Z(1/2+it, 1/2+iu; \psi, \psi')|^2 du \, dt \ll (Y_1Y_2)^{1+\varepsilon}
\end{displaymath}
for any $Y_1, Y_2 \geq 1$. 
\end{abstract}


\maketitle

\section{Introduction}

Subconvexity bounds for $L$-functions is one of the central topics in analytic number theory with deep and sometimes unexpected applications in almost every branch of number theory.  H.\ Weyl, almost a century ago, was the first to prove a subconvex bound for the Riemann zeta-function: $\zeta(1/2 + it) \ll |t|^{1/6 +\varepsilon}$. Since then,  powerful methods from various areas have been developed culminating in a complete solution of the subconvexity problem for $L$-functions on $GL_1$ and $GL_2$ due to Michel and Venkatesh \cite{MV}. One would hope that new methods emerge that will enable subconvex bounds for more general $L$-functions. This may include automorphic $L$-functions of higher rank groups (some deep, but sporadic results are already available, e.g.\ \cite{HM}, \cite{Li}, \cite{Ve}), $L$-functions without Euler product (e.g.\ $L$-functions of half-integral weight modular forms) or multiple $L$-functions, that is, $L$-functions whose coefficients are again $L$-functions. Multiple $L$-functions have become more and more frequent in analytic number theory, and have quite recently proved to be a very powerful and elegant tool that in some cases can prove results that are not (yet) available with other methods, see e.g.\ \cite{DGH}. A good understanding of the more subtle anayltic properties of multiple Dirichlet series   would be very desirable. The question for subconvexity bounds for double Dirichlet series has been raised explicitly in \cite{HK} in connection with non-vanishing results for quadratic twists.   

The aim of this paper is to show the first subconvex bound for a multiple Dirichlet series in a relatively simple situation. For $\Re s$ and $\Re w$ sufficiently large and two real characters $\psi, \psi'$ of conductor dividing 8 we define
\begin{equation}\label{defZ}
  Z(s, w; \psi, \psi') := \zeta_2(2s+2w-1) \sum_{d \text{ odd}} \frac{L_2(s, \chi_d\psi) \psi'(d)}{d^w}
\end{equation}
where $\chi_d = \left(\frac{d}{.}\right)$ and here and henceforth the subscripts 2 denote that the Euler factors at 2 have been removed. This type of series was first considered in \cite{GH}, although not from the point of view of double Dirichlet series. It has two more or less obvious functional equations: the functional equation for $L(s, \chi)$ yields a functional equation sending $s$ to $1-s$, and interchanging the order of summation yields a functional equation interchanging $s$ and $w$. These two functional equations generate the dihedral group $D_{12}$ of order 12, and successive application of the functional equations yields the meromorphic continuation of $Z(s, w, \psi, \psi')$ as a function of two complex variables with polar lines at most at $s=1$, $w=1$ and $s+w = 3/2$. It is a priori not completely obvious what the convexity bound in this situation is, because running a convexity argument for $Z(s, w; \psi, \psi')$ depends on what we assume on the coefficients $L(s, \chi)$ in the region of absolute convergence of $Z(s, w; \psi, \psi')$. If we assume the convexity bound for $L(s, \chi_d)$ in $s$ together with the Lindel\"of hypothesis on average over $d$, that is,
\begin{displaymath}
  \sum_{d \leq X} |L(s, \chi_d)| \ll (X|s|^{1/4})^{1+\varepsilon}, \quad \Re s =1/2,
\end{displaymath}
(cf.\ \eqref{HBest} below), the convexity bound for $Z(s, w; \psi, \psi')$ is 
\begin{equation}\label{conv}
  Z(s, w; \psi, \psi') \ll |sw(s+w)|^{1/4+\varepsilon}
\end{equation}
for $\Re s= \Re w = 1/2$. This is indeed a natural candidate for the convexity bound, since 
\begin{displaymath}
 \Gamma\left(\frac{s}{2}\right)\Gamma\left(\frac{s+w-1/2}{2}\right)\Gamma\left(\frac{w}{2}\right) Z(s, w; \psi, \psi')
\end{displaymath} 
  is roughly invariant under $(s, w) \mapsto (1-s, 1-w)$, see \eqref{feq3} below, hence it is reasonable to define the "analytic conductor" of $Z(1/2 + it, 1/2+iu; \psi, \psi')$ by
\begin{equation}\label{anacon}
  C(u, t) :=  \left|\frac{1}{4} + \frac{it}{2} \right|\cdot  \left|\frac{1}{4} + \frac{i(u+t)}{2} \right| \cdot \left|\frac{1}{4} + \frac{iu}{2} \right|.
\end{equation}
We shall prove the following uniform subconvexity bound.


\begin{satz} One has
\begin{displaymath}
  Z(s, w, \psi, \psi') \ll |sw(s+w)|^{1/6+\varepsilon}
\end{displaymath}
for $\Re s = \Re w = 1/2$ any $\varepsilon > 0$.
\end{satz}

This matches the quality of Weyl's bound for the Riemann zeta-function and the corresponding estimate for $L$-functions attached to modular forms on $GL_2$ due to A.\ Good \cite{Go}. Theorem 1 is the first subconvex  bound for a multiple Dirichlet series, and it seems to be the first subconvex result in the literature for an $L$-series that is not a (linear combination of) $L$-series with Euler product. 

The bound of Theorem 1 is non-trivial even in a one-variable situation.  Specializing to $s = 1/2$,   one gets an ordinary Dirichlet series (without an Euler product) with coefficients given by central $L$-values. Then $\Gamma(w)Z(1/2, w; \psi, \psi')$ is roughly invariant under $w \mapsto 1-w$, hence the standard convexity bound in one variable is $|w|^{1/2+\varepsilon}$ on the critical line which coincides with \eqref{conv} in this case.  Theorem 1 implies the  subconvex bound $Z(1/2, w; \psi, \psi')\ll |w|^{1/3+\varepsilon}$. 

Another interesting case comes from the specialization is $s+w =1$, that is, $s = 1/2 + it$, $w = 1/2 - it$. Of course, $Z(s, 1-s; \psi, \psi')$ exists only by analytic continuation. This is a situation where  the analytic conductor  \eqref{anacon} is unusually small due to a special configuration. This is a well-known phenomenon that occurs for instance with $L$-functions of Maa{\ss } forms with spectral parameter $t$ in the neighbourhood of  the point $1/2 \pm it$. Such effects have quite interesting consequences, see for example \cite{Y} or \cite{Bl}. For $GL(2)$-Maa{\ss} forms, the subconvexity problem in this situation has been solved only recently in \cite{MV}. Theorem 1 above yields $Z(1/2 + it, 1/2 - it; \psi, \psi') \ll (1+|t|)^{1/3+\varepsilon}$ while the convexity bound is $(1+|t|)^{1/2+\varepsilon}$. 

One may speculate if a Lindel\"of type bound holds for $Z(s, w; \psi, \psi')$. In this direction we prove

\begin{satz} For $Y_1, Y_2 \geq 1$ one  has
\begin{displaymath}
  \int_{-Y_1}^{Y_1} \int_{-Y_2}^{Y_2} |Z(1/2 + it, 1/2+iu; \psi, \psi')|^2 du \, dt\ll (Y_1Y_2)^{1+\varepsilon}.
\end{displaymath}
for any $\varepsilon > 0$. 
\end{satz}




Both Theorem 1 and 2 can be extended to Dirichlet series $Z(s, w; \psi, \psi')$ where $\psi, \psi'$ are allowed to have any (fixed) ramification. \\

The proofs of Theorems 1 and 2 start by writing down approximate functional equations for $Z(s, w; \psi, \psi')$. For multiple Dirichlet series we have several choices. Since $\Gamma(s/2)\Gamma(w/2)\Gamma((s+w)/2) Z(s, w; \psi, \psi')$ is roughly invariant under $(s, w) \mapsto (1-s, 1-w)$, see \eqref{feq3} below, one has the simplest approximate functional equation of the type
\begin{equation}\label{test}
  Z(1/2 + it, 1/2+iu;\psi, \psi') \approx \sum_{d \leq (u(u+t))^{1/2}} \frac{L_2(1/2+it, \chi_d\psi) \psi'(d)}{d^{1/2+iu}},
\end{equation}
cf.\ Lemma 3.  An average bound for $L$-values (see \eqref{HBest}) then recovers, as usual,  the convexity bound.  
We can now insert another approximate functional into the numerator $L_2(1/2+it, \chi_d \psi)$ getting something roughly of the form
\begin{equation}\label{heu1}
  Z(1/2+it, 1/2+iu;\psi, \psi')  \ll \sum_{d \sim P} \sum_{n \sim Q} \frac{\chi_d(n)}{n^{1/2+it} d^{1/2+iu}},  \quad P \leq (u(u+t))^{1/2},  \quad Q \leq \sqrt{Pt}. 
\end{equation} 
This gives good bounds  if $P$ happens to be small. 
 For large $P$ we can successfully apply Poisson-summation in the long $d$-variable (by quadratic reciprocity it is self-dual at about $\sqrt{Qu}$) and find something roughly of the form 
 \begin{equation}\label{heu2}
  Z(1/2+it, 1/2+iu;\psi, \psi')  \ll \sum_{d \sim Qu/P} \sum_{n \sim Q} \frac{\chi_d(n)}{n^{1/2+i(t+u)} d^{1/2-iu}}.
\end{equation} 
Theorem 1 follows now  from Heath-Brown's large sieve estimate \eqref{HBbilinear} which allows to bound efficiently bilinear sums in $\chi_d(n)$.  
Theorem 2 follows from \eqref{heu1} and \eqref{heu2} together with standard bounds for Dirichlet polynomials. 
The above approach   based on Poisson summation serves as a good heuristic, but  has to be modified. Not all numbers are squarefree and $\equiv 1 $ (mod 4), and a rigorous argument along these lines would face similar substantial difficulties as in \cite{HB}. However, Poisson summation in the $d$-variable can be mimiced by applying a suitable functional equation of $Z(s, w; \psi, \psi')$ sending $(s, w)$ to $(s+w-1/2, 1-w)$. Lemma 4 and 5 will provide the crucial bounds that correspond roughly to \eqref{heu1} and \eqref{heu2}. They could be turned into equations with small error term (what one might call an approximate functional equation), but the formulas would become even more cumbersome. We remark that it is important for the proofs of Theorems 1 and 2 that the variables are cleanly separated which makes the argument a little more technical than the  heuristic approach.   \\

\textbf{Notation.} Most of the notation is standard. We recall that $\varepsilon$ denotes a sufficiently small positive constant, not necessarily the same at each occurrence. Similarly, $A$ denotes a sufficiently large positive constant, not necessarily the same at each occurrence. The notation $x \sim X$ means $X \leq x \leq 2X$. All implied constants may depend on $\varepsilon$ and/or $A$ even if not explicitly specified.

\section{Preparatory material}

\subsection{Real characters}

We follow the notation of \cite{DGH}. Let $d$ and $n$ be odd positive integers that we decompose uniquely as $d = d_0d_1^2$ with $\mu^2(d_0) = 1$ and $n = n_0n_1^2$ with $\mu^2(n_0) = 1$. We write
\begin{equation}\label{change}
  \chi_d(n) = \left(\frac{d}{n}\right) = \tilde{\chi}_n(d).
\end{equation}
The character $\chi_d$ is the Jacobi-Kronecker symbol of conductor $d_0$ if $d \equiv 1$ (mod 4) and $4d_0$ if $d \equiv 3$ (mod 4). We have
\begin{equation}\label{2factor}
  \chi_d(2) = \begin{cases} 1, & d \equiv 1 \, (\text{mod } 8);\\
  -1, & d \equiv 5 \, (\text{mod } 8);\\
  0, & d \equiv 3 \, (\text{mod } 4).
  \end{cases}
\end{equation}
and $\chi_d(-1) = 1$, that is, $\chi_d$ is even. By quadratic reciprocity we have
\begin{equation}\label{quadrrec}
 \tilde{\chi}_n = \begin{cases}
   \chi_n, &\quad  n \equiv 1 \, (\text{mod  }4);\\
   \chi_{-n}, & \quad n \equiv 3 \, (\text{mod  }4).
   \end{cases}
\end{equation}     

Let $\psi_1, \psi_{-1}, \psi_2, \psi_{-2}$ denote the four characters modulo 8, that is, $\psi_1$ is the trivial character, $\psi_{-1}$ is induced from the non-trivial character modulo 4, $\psi_2(n) = 1$ if and only if $n \equiv 3$ or $5$ (mod 8) and $\psi_{-2}(n) = 1$ if and only if $n \equiv 5$ or $7$ (mod 8). 

By considering $\chi_{d_0} \psi$ for  odd squarefree $d$ and $\psi$ (mod 8) we can construct all real primitive characters. The $L$-series satisfies a functional equation
\begin{equation}\label{functeq}
  L(s, \chi_{d_0} \psi)  =  \left(\frac{\delta_0}{\pi}\right)^{\frac{1}{2} - s} \frac{\Gamma(\frac{1-s+\kappa}{2})}{\Gamma(\frac{s+\kappa}{2})} L(1-s, \chi_{d_0} \psi)
\end{equation}
where
\begin{equation}\label{fepara}
  \kappa = \begin{cases}0, & \psi = \psi_1 \text{ or } \psi_2,\\ 1,& \psi = \psi_{-1} \text{ or } \psi_{-2}; \end{cases} \quad \delta_0 = \begin{cases} d_0, & \psi = \psi_1, d \equiv 1 \, (4) \text{ or } \psi = \psi_{-1}, d \equiv 3 \, (4),\\ 4d_0, & \psi = \psi_1, d \equiv 1 \, (4) \text{ or } \psi = \psi_{-1}, d \equiv 3 \, (4),
  \\ 8d_0, &  \psi = \psi_2 \text{ or } \psi_{-2}.
  \end{cases}
\end{equation}
This gives an approximate functional equation \cite[Theorem 5.3]{IK}
\begin{equation}\label{approxL}
  L(1/2+it, \chi_{d_0} \psi) = \sum_{n} \frac{(\chi_{d_0}\psi)(n)}{n^{1/2+it}} G^{(\psi)}_t\left(\frac{n}{\sqrt{\delta_0}}\right) + \lambda(t, \delta_0)\sum \frac{(\chi_{d_0}\psi)(n)}{n^{1/2-it}} G^{(\psi)}_{-t}\left(\frac{n}{\sqrt{\delta_0}}\right) 
\end{equation}
where $|\lambda(t, \delta_0)|=1$ and, for an arbitrary number $A > 5$, 
\begin{equation}\label{defG}
  G^{(\psi)}_t(\xi) = \frac{1}{2\pi i} \int_{(2)}\left(\cos\frac{\pi s}{4A}\right)^{-4A} \frac{\Gamma(\frac{\frac{1}{2} + it +s +\kappa}{2})}{\Gamma(\frac{\frac{1}{2} + it +\kappa}{2})} \xi^{-s} \frac{ds}{s} \ll \Biggr(1+\frac{\xi}{\sqrt{1+|t|}}\Biggl)^{-A},
\end{equation}
cf.\ \cite[Proposition 5.4]{IK}. Most of the time we shall deal with non-fundamental discriminants, therefore we note that 
\begin{equation}\label{nonfund}
  L_2(1/2 + it, \chi_d\psi) = \prod_{p \mid 2d_1} \left(1-\frac{(\chi_{d_0}\psi)(p)}{p^{1/2+it}}\right) L(1/2 + it, \chi_{d_0}\psi).
\end{equation}

The paper rests crucially on Heath-Brown's large sieve estimate for quadratic characters.  Here we need the following two corollaries: Let $(a_m), (b_n)$ be sequences of complex numbers numbers of absolute value at most 1, then \cite[Corollary 4]{HB} states
\begin{displaymath}
  \sum_{\substack{m\leq M\\ m \text{ odd}}} \sum_{n \leq N} a_mb_n \left(\frac{n}{m}\right) \ll_{\varepsilon} (MN)^{\varepsilon}(MN^{1/2} + M^{1/2}N)
\end{displaymath}
for any $\varepsilon >0$. We will use it in the following form: if $\tilde{a}_m, \tilde{b}_m \ll m^{-1/2+\varepsilon}$, then 
\begin{equation}\label{HBbilinear}
  \sum_{\substack{m\leq M\\ m \text{ odd}}} \sum_{n \leq N} \tilde{a}_m\tilde{b}_n \left(\frac{n}{m}\right) \ll_{\varepsilon} (M+N)^{1/2+\varepsilon}.
\end{equation}
This follows directly from the preceding inequality after cutting into dyadic pieces.  Secondly, \cite[Theorem 2]{HB} states
\begin{displaymath}
  \sum_{\substack{d_0 \leq X\\ d_0 \text{ odd, squarefree}}} |L(s, \chi_{d_0}\psi)|^4 \ll (X(|s|+1))^{1+\varepsilon}, \quad \sigma \geq 1/2,
\end{displaymath}
from which one obtains by \eqref{nonfund} and H\"older's inequality the bound
\begin{equation}\label{HBest}
  \sum_{\substack{d \leq X\\ d  \text{ odd}}} |L_2(s, \chi_{d}\psi)| \ll (X|s|^{1/4})^{1+\varepsilon}, \quad \sigma \geq 1/2;
\end{equation}
recall that the subscript $2$ indicates the removal of the Euler factor at 2. Heath-Brown's original  bound is somewhat stronger, and for the purpose of this paper a second moment would suffice. 


\subsection{Special functions}

We recall Stirling's formula in the following form:  For $s, z \in \Bbb{C}$ with $\Re(s+z) \geq 1/10$ we have the uniform bound
\begin{equation}\label{stir}
  \frac{\Gamma(s+z)}{\Gamma(\bar{s}-z)} \ll_{\Re s, \Re z}  |s+z|^{2\Re z} 
\end{equation}
and \cite[p.\ 100]{IK}
\begin{equation}\label{stir2}
  \frac{\Gamma(s+z)}{\Gamma(s)} \ll_{\Re s, \Re z}  (1+|s|)^{\Re z} \exp\left(\frac{\pi}{2} |z|\right).  
\end{equation}
For future reference we remark that
\begin{equation}\label{cot}
\frac{\Gamma(\frac{2-z}{2})}{\Gamma(\frac{z+1}{2})} = \frac{\Gamma(\frac{1-z }{2})}{\Gamma(\frac{z }{2})} \cot\left(\frac{\pi z}{2}\right), \quad z \in \Bbb{C}. 
\end{equation}
Away from poles, we have the uniform asymptotic formula 
\begin{equation}\label{cotbound}
   \cot(x+iy) = -i \, \text{sign}(y)  + O(e^{-2|y|}), \quad \min_{k \in \Bbb{Z}}|z-\pi k|  \geq 1/10. 
\end{equation}

\subsection{Multiple Dirichlet polynomials} For the proof of Theorem 2 we will need the following lemma.

\begin{lemma} Let $D, N > 0$, $Y_1, Y_2 \geq 1$, $X := Y_1Y_2DN$, $\varepsilon > 0$.   Let $W_1, W_2$ be two fixed smooth functions with support in $[-2, 2]$. For $n, d \in \Bbb{N}$ let $f(d, n)$ be any complex numbers with absolute value at most 1. Then
\begin{equation}\label{diri1}
\int\int W_1\left(\frac{t}{Y_1}\right) W_2\left(\frac{u}{Y_2}\right)\Biggl|\sum_{\substack{d \sim D\\n \sim N}} \frac{f(d, n)}{n^{it}d^{iu}}\Biggr|^2 du\, dt \ll  X^{\varepsilon} DNY_1Y_2 \left(1+\frac{N}{Y_1}\right)\left(1+\frac{D}{Y_2}\right)
\end{equation}
and
\begin{equation}\label{diri2}
\begin{split}
& \int\int W_1\left(\frac{t}{Y_1}\right) W_2\left(\frac{u}{Y_2}\right)\Biggl|\sum_{\substack{d \sim D\\n \sim N}} \frac{f(d, n)}{n^{iu}d^{-i(u+t)}}\Biggr|^2 du\, dt \\
& \ll  X^{\varepsilon}\left(NDY_1Y_2 +ND^2\min(Y_1, Y_2) + N^2DY_1 + (ND)^2\right). 
\end{split}
\end{equation}
\end{lemma}

\textbf{Proof.}  Without loss of generality we can assume $D, N \geq 1$, otherwise the $d, $-sums vanish. Opening the square, we bound the left hand side of \eqref{diri1} by
\begin{displaymath}
  \sum_{\substack{d_1, d_2 \sim D \\ n_1, n_2 \sim N}} \left| \int\int W_1\left(\frac{t}{Y_1}\right) W_2\left(\frac{u}{Y_2}\right) \left(\frac{n_2}{n_1}\right)^{it} \left(\frac{d_2}{d_1}\right)^{iu} du\, dt\right|. 
\end{displaymath}
Integrating by parts sufficiently often, we can assume, up to an error $O(X^{-A})$,  
\begin{displaymath}
  n_2 = n_1 \left(1+O(X^{\varepsilon}Y_1^{-1})\right), \quad d_2 = d_1 \left(1+O(X^{\varepsilon}Y_2^{-1})\right),
\end{displaymath}
and \eqref{diri1} follows immediately. The second part requires a slightly more careful argument. Again we bound the left hand side of \eqref{diri2} by
\begin{displaymath}
  \sum_{\substack{d_1, d_2 \sim D \\ n_1, n_2 \sim N}} \left| \int\int W_1\left(\frac{t}{Y_1}\right) W_2\left(\frac{u}{Y_2}\right) \left(\frac{n_2d_1}{n_1d_2}\right)^{iu} \left(\frac{d_1}{d_2}\right)^{it} du\, dt\right|, 
\end{displaymath}
and we can restrict the summation to
\begin{displaymath}
  d_1 = d_2 \left(1+O\left(X^{\varepsilon} Y_1^{ -1}\right)\right), \quad n_2d_1 = n_1d_2 \left(1+O\left(X^{\varepsilon}Y_2^{ -1}\right)\right). 
\end{displaymath} 
Hence the left hand side of \eqref{diri2} is at most 
\begin{equation}\label{count0}
  Y_1Y_2   \#\mathcal{A}
   \end{equation}
  where $\mathcal{A}$ is the set of all 6-tuples $(d_1, d_2, n_1, n_2, a, b)\in \Bbb{Z}^6$ satisfying
  \begin{align}\label{conditions}
\nonumber &  d_1, d_2 \sim D,\quad n_1, n_2 \sim N, \quad |a| \ll D X^{\varepsilon} Y_1^{-1} , \quad |b| \ll NDX^{\varepsilon}Y_2^{-1},\\ 
& d_1 = d_2 + a, \quad n_2d_1 = n_1d_2 + b.
  \end{align}
The number of such 6-tuples with $n_1 = n_2$ is 
\begin{equation}\label{count1}
  \ll X^{\varepsilon} N D \left(1+ D \max(Y_1, Y_2)^{-1}\right).
\end{equation}   
Let us now assume $n_1 \not=n_2$. We substitute the first equation in \eqref{conditions} into the second and write $n_3 = n_1-n_2 \not=0$. Hence $\#\mathcal{A}$ is at most the number of 5-tuples $(d_1, n_1, n_3, a, b)$ satisfying  
  \begin{displaymath}
  \begin{split}
&  d_1, d_2 \sim D,\quad n_1 \sim N, \quad  0 < |n_3| \leq N, \quad |a| \ll DX^{\varepsilon} Y_1^{-1} , \quad |b| \ll NDX^{\varepsilon}Y_2^{-1},  \quad n_3d_1 = n_1a + b 
  \end{split}
  \end{displaymath}  
which, by a divisor argument, is at most 
\begin{equation}\label{count2}
 \ll X^{\varepsilon} N \left(1+\frac{D}{Y_1}\right) \left(1+\frac{ND}{Y_2}\right). 
\end{equation}  
We substitute \eqref{count1} and \eqref{count2} into \eqref{count0} and arrive at the right hand side of \eqref{diri2}. 

\section{Functional equation and meromorphic continuation}
 
The aim of this section is establish the meromorphic continuation and the functional equations of the double Dirichlet series $Z(s, w; \psi, \psi')$ defined in \eqref{defZ}. We will treat these 16 series simultaneously and introduce the following notation: Let 
\begin{displaymath}
\textbf{Z}(s, w; \psi) = \left(\begin{array}{l} Z(s, w; \psi, \psi_1)\\ Z(s, w; \psi, \psi_{-1})\\ Z(s, w; \psi, \psi_2)\\ Z(s, w; \psi, \psi_{-2}) \end{array}\right), \quad \textbf{Z}(s, w) := \left(\begin{array}{l} \textbf{Z}(s, w; \psi_1)\\ \textbf{Z}(s, w; \psi_{-1})\\ \textbf{Z}(s, w;  \psi_2)\\ \textbf{Z}(s, w;  \psi_{-2}) \end{array}\right),
\end{displaymath}
so $\textbf{Z}(s, w)$ is a column vector with 16 entries. We have the following lemma.
 
\begin{lemma} The functions $(s-1)(w-1)(s+w-3/2)Z(s, w; \psi, \psi')$  can be extended holomorphically  to all of $\Bbb{C}^2$.  They are of at most  polynomial growth in $\Im s$ and $\Im w$ in the sense that for any $C_1 > 0$ there is a constant $C_2>0$ such that $(s-1)(w-1)(s+w-3/2)Z(s, w; \psi, \psi') \ll ((1+\Im s)(1+\Im w))^{C_2}$ whenever $|\Re s|, |\Re w| \leq C_1$.  Moreover, there are 16-by-16 matrices $\mathcal{A}$ and $\mathcal{B}(s)$ given by \eqref{defA} and \eqref{defB} below, such that 
\begin{equation}\label{feqs1}
\textbf{Z}(s, w) = \mathcal{A} \, \textbf{Z}(w, s)
\end{equation}
and
\begin{equation}\label{feqs2}
   \textbf{Z}(s, w) = \mathcal{B}(s) \textbf{Z}(1-s, s+w -1 /2).
\end{equation}
\end{lemma} 
 
\textbf{Proof.} This is essentially known and follows the procedure outlined in \cite[Section 4]{DGH}. For convenience, we give the complete argument and provide explicit formulas.  

We start with following two expressions for $Z(s, w; \psi, \psi')$, initially valid for $\Re s, \Re w$ sufficiently large. On the one hand we have 
\begin{equation}\label{fe2a}
\begin{split}
Z(s, w;  \psi, \psi') & = \zeta_2(2s+2w-1) \sum_{\substack{d_0 \text{ odd}\\ \mu^2(d_0) = 1}} \frac{L_2(s, \chi_{d_0}\psi) \psi'(d_0)}{d_0^w} \sum_{d_1 \text{ odd}} \frac{1}{d_1^{2w}}\prod_{p \mid d_1} \left(1-\frac{(\chi_{d_0}\psi)(p)}{p^s}\right)\\
& = \zeta_2(2s+2w-1) \sum_{\substack{d_0 \text{ odd}\\ \mu^2(d_0) = 1}} \frac{L_2(s, \chi_{d_0}\psi) \psi'(d_0) \zeta_2(2w)}{d_0^w L_2(s + 2w, \chi_{d_0} \psi)}.
\end{split}
\end{equation}
The right-hand side of  \eqref{fe2a} is, by \eqref{HBest} together with \eqref{functeq} for $\Re s < 1/2$, absolutely and locally uniformly convergent and hence holomorphic  in 
\begin{displaymath}
  R_1 :=  \{(s, w) \mid \Re w > 1\} \cap \{(s, w) \mid \Re s + \Re w > 3/2\}
\end{displaymath}
with the exception of a polar line at $s=1$ if $\psi = \psi_1$ is trivial, and it is of moderate growth in $\Im s$, $\Im w$ in this region. 

On the other hand we have
 \begin{equation}\label{feq1}
\begin{split}
  Z(s, w; & \psi, \psi') = \zeta_2(2s+2w-1) \sum_{d \text{ odd}} \frac{L_2(s, \chi_d\psi) \psi'(d)}{d^w}\\
  & = \zeta_2(2s+2w-1) \sum_{d, n \text{ odd}} \frac{\chi_d(n) \psi(n) \psi'(d)}{d^wn^s} = \zeta_2(2s+2w-1)\sum_{n \text{ odd}} \frac{L_2(w, \tilde{\chi}_n \psi') \psi(n)}{n^s}
\end{split}  
\end{equation} 
where we used \eqref{change}. The two equalities \eqref{fe2a} and \eqref{feq1} together with \eqref{2factor} - \eqref{fepara} yield now readily the two matrices $\mathcal{A}$ and $\mathcal{B}(s)$. One way to construct the matrices explicitly  is as follows. For a character $\psi$ mod 8 and a residue class $\eta \in \{1, 3, 5, 7\}$ mod 8 let
\begin{equation}\label{defY}
  Y_{\eta}(s, w, \psi) =  \zeta_2(2s+2w-1) \sum_{d \equiv \eta \, (8)} \frac{L_2(s, \chi_d\psi) \psi'(d)}{d^w} = \frac{1}{4} \sum_{\psi'} \psi'(\eta) Z(s, w, \psi, \psi') 
\end{equation}
and 
\begin{displaymath}
\textbf{Y}(s, w; \psi) = \left(\begin{array}{l} Y(s, w; \psi, 1)\\ Y(s, w; \psi, 3)\\ Y(s, w; \psi, 5)\\ Y(s, w; \psi, 7) \end{array}\right),  \quad \textbf{Y}(s, w) := \left(\begin{array}{l} \textbf{Y}(s, w; \psi_1)\\ \textbf{Y}(s, w; \psi_{-1})\\ \textbf{Y}(s, w;  \psi_2)\\ \textbf{Y}(s, w;  \psi_{-2}) \end{array}\right),
\end{displaymath}
Moreover, let $\tilde{\textbf{Y}}(s, w)$ and $\tilde{\textbf{Z}}(s, w)$ the same 16-by-16 vectors as $\textbf{Y}(s, w)$ and $\textbf{Z}(s, w)$, resp.\, except that  in numerator of each component of $\tilde{\textbf{Y}}(s, w)$ and $\tilde{\textbf{Z}}(s, w)$   the character $\chi_d$ in $L_2(s, \chi_d\psi)$ is replaced with $\tilde{\chi}_d$. Then \eqref{defY} gives readily a relation $\textbf{Z}(s, w) = \mathcal{M}_1\textbf{Y}(s, w)$ for a 16-by-16 matrix $\mathcal{M}_1$ consisting of four identical 4-by-4 blocks on the diagonal. Next, by \eqref{quadrrec} we find a matrix $\mathcal{M}_2$ such that    $\textbf{Y}(s, w) = \mathcal{M}_2 \tilde{\textbf{Y}}(s, w)$. Now we use the equation \eqref{feq1} to get a functional equation $\textbf{Z}(s, w) = \mathcal{M}_3 \tilde{\textbf{Z}}(w, s)$. Finally, applying the functional equation \eqref{functeq} and \eqref{fepara} together with \eqref{2factor} to \eqref{fe2a}, we find a diagonal   matrix $\mathcal{M}_4(s)$ such that $\textbf{Y}(s, w) = \mathcal{M}_4(s) \textbf{Y}(1-s, s+w -1/2)$. Note that the map  $(s, w) \mapsto (1-s, s+w-1/2)$  leaves $s+2w$ invariant and interchanges $2s+2w-1$ and $2w$. Putting together these four matrix equations, we get \eqref{feqs1} with $\mathcal{A} = \mathcal{M}_3\mathcal{M}_1\mathcal{M}_2\mathcal{M}_1^{-1}$ and \eqref{feqs2} with  $\mathcal{B}(s) = \mathcal{M}_1\mathcal{M}_4(s) \mathcal{M}_1^{-1}$, explicitly\begin{equation}\label{defA}
  \mathcal{A} = \frac{1}{2} \left(\begin{smallmatrix}1 & 1 & 0 & 0 & 1 & -1 & 0 & 0 & 0 & 0 & 0 & 0 & 0 & 0 & 0 & 0\\1 & -1 & 0 & 0 & 1 & 1 & 0 & 0 & 0 & 0 & 0 & 0 & 0 & 0 & 0 & 0\\ 0 & 0 & 0 & 0 & 0 & 0 & 0 & 0 &  1 & 1 & 0 & 0 & 1 & -1 & 0 & 0\\ 0 & 0 & 0 & 0 & 0 & 0 & 0 & 0 &  1& -1 & 0 & 0 &1 &  1 & 0 & 0\\ 1 & 1 & 0 & 0 & -1 & 1 & 0 & 0 & 0 & 0 & 0 & 0 & 0 & 0 & 0 & 0\\ -1 & 1 & 0 & 0 & 1 & 1 & 0 & 0 & 0 & 0 & 0 & 0 & 0 & 0 & 0 & 0\\ 0 & 0 & 0 & 0 & 0 & 0 & 0 & 0 &  1 & 1 & 0 & 0 & -1 & 1 & 0 & 0\\ 0 & 0 & 0 & 0 & 0 & 0 & 0 & 0 &  -1& 1 & 0 & 0 &1 &  1 & 0 & 0\\
0 & 0 &  1 & 1 & 0 & 0 & 1 & -1 & 0 & 0 & 0 & 0 & 0 & 0 & 0 & 0  \\0 & 0 & 1 & -1 & 0 & 0 & 1 & 1 & 0 & 0 & 0 & 0 & 0 & 0 & 0 & 0 \\ 0 & 0 & 0 & 0 & 0 & 0 & 0 & 0 & 0 & 0 &  1 & 1 & 0 & 0 & 1 & -1  \\0 & 0 &  0 & 0 & 0 & 0 & 0 & 0 & 0 & 0 &  1& -1 & 0 & 0 &1 &  1  \\ 0 & 0 & 1 & 1 & 0 & 0 & -1 & 1 & 0 & 0 & 0 & 0 & 0 & 0 & 0 & 0 \\0 & 0 &  -1 & 1 & 0 & 0 & 1 & 1 & 0 & 0 & 0 & 0 & 0 & 0 & 0 & 0 \\ 0 & 0 & 0 & 0 & 0 & 0 & 0 & 0 & 0 & 0 &  1 & 1 & 0 & 0 & -1 & 1 \\ 0 & 0 & 0 & 0 & 0 & 0 & 0 & 0 & 0 & 0 &  -1& 1 & 0 & 0 &1 &  1 
  \end{smallmatrix}\right) \in \Bbb{C}^{16\times 16}.
\end{equation}
and   
 \begin{equation}\label{defB}
 \begin{split}
&  \mathcal{B}(s) = \left(\begin{matrix} \mathcal{B}_1(s) & & & \\ & \mathcal{B}_2(s) & & \\ & & \mathcal{B}_3(s) & \\ & & & \mathcal{B}_4(s)\end{matrix}\right),\\
 &\mathcal{B}_1(s) = \frac{\pi^{s-\frac{1}{2}} \Gamma(\frac{1-s}{2})}{(4^s-4) \Gamma(\frac{s}{2})} \left(\begin{smallmatrix} -4^{1-s} & 4^{1-s}-2 & 2^{1-s}-2^s & 2^{1-s} - 2^s\\  4^{1-s}-2 & - 4^{1-s} & 2^{1-s}-2^s & 2^{1-s} - 2^s \\ 2^{1-s}-2^s & 2^{1-s} - 2^s & -4^{1-s} & 4^{1-s}-2 \\ 2^{1-s}-2^s & 2^{1-s} - 2^s & 4^{1-s}-2&  -4^{1-s}  \end{smallmatrix}\right),\\
 & \mathcal{B}_2(s) = \frac{\pi^{s-\frac{1}{2}} \Gamma(\frac{2-s}{2})}{(4^s-4) \Gamma(\frac{s+1}{2})} \left(\begin{smallmatrix} -4^{1-s} & 2-4^{1-s} & 2^{1-s}-2^s & 2^s - 2^{1-s} \\  2-4^{1-s} & - 4^{1-s} & 2^s - 2^{1-s} & 2^{1-s} - 2^s \\ 2^{1-s}-2^s & 2^s - 2^{1-s}  & -4^{1-s} & 2- 4^{1-s} \\ 2^s - 2^{1-s} & 2^{1-s} - 2^s & 2-4^{1-s}&  -4^{1-s}  \end{smallmatrix}\right), \\ 
& \mathcal{B}_3(s) = \left(\frac{\pi}{8}\right)^{s-\frac{1}{2}}\frac{\Gamma(\frac{1-s}{2})}{\Gamma(\frac{s}{2})} I_4, \quad \mathcal{B}_4(s) = \left(\frac{\pi}{8}\right)^{s-\frac{1}{2}}\frac{\Gamma(\frac{2-s}{2})}{\Gamma(\frac{s+1}{2})} I_4.
\end{split}
\end{equation}
The two functional equations \eqref{feqs1} and \eqref{feqs2} are involutions and generate the dihedral group of order 12.   The exact shape of the matrices $\mathcal{A}$ and $\mathcal{B}(s)$ is not   important, but we note that the entries of $\mathcal{B}(s)$ are 
\begin{equation}\label{propB}
  \text{holomorphic and of moderate growth in }\Im s \text{ if } \Re s < 1.
\end{equation}   
and 
\begin{equation}\label{propB2}
  \mathcal{B}_1(0) = 0.
\end{equation}

We proceed to continue $Z(s, w; \psi, \psi')$ meromorphically. Let $\alpha(s, w) = (w, s)$ and $\beta(s, w) = (1-s, s+w-1/2)$. Since $\alpha(R_1) \cap R_1$ is an open set in $\Bbb{C}^2$, we can apply \eqref{feqs1} to continue $Z(s, w; \psi, \psi')$ to the region
\begin{displaymath}
  R_2 := \alpha(R_1) \cup R_1 = \{(s, w) \mid \Re s+\Re w > 3/2, \max(\Re s, \Re w) > 1\}.  
\end{displaymath}
with moderate growth in $\Im w$ and $\Im s$ and polar lines at most at $s=1$ and $w=1$. Next, since $\beta(R_2) \cap R_2$ is open in $\Bbb{C}^2$, we can apply \eqref{feqs2} and continue $Z(s, w;\psi, \psi')$ to 
\begin{displaymath}
  R_3 := \beta(R_2) \cup R_2 = R_2 \cup \{(s, w) \mid \Re s < 0, \Re w > 1\}
\end{displaymath}
By \eqref{propB}, $Z(s, w; \psi, \psi')$ is of moderate growth in $R_3$, and the only possible singularities in $R_3 \setminus R_2$ can occur at $\beta\{(1, w) \mid w \in \Bbb{C}\} = \{(0, w) \mid w \in \Bbb{C}\}$ and  $\beta\{(s, 1) \mid s \in \Bbb{C}\} = \{(s, w) \mid s + w = 3/2\}$. By \eqref{propB2}, the first case cannot occur. Next we apply \eqref{feqs1} again getting a continuation to 
\begin{displaymath}
  R_4 := \{(s, w) \mid \max(\Re s, \Re w) > 1, \Re s + \Re w > 3/2 \text{ if } \Re s, \Re w \geq 0\}. 
\end{displaymath}
Finally, we apply once again \eqref{feqs2} (getting no new singularities since the line $s+w= 3/2$ is mapped to $w=1$) and \eqref{feqs1}. In this way we establish the meromorphic continuation with moderate growth to all of $\Bbb{C}^2$ with the exception of the tube 
\begin{displaymath}
  R^{\ast}  := \{(s, w) \mid (\Re s, \Re w) \in \Omega\} \subseteq \{(s, w) : |\Re s|^2 + |\Re w|^2 \leq 3\}
\end{displaymath}  
where $\Omega \subseteq \Bbb{R}^2$ is the closed 12-gon with vertices
\begin{displaymath}
 (1, 1), (1/2, 1), (0, 3/2), (0, 1), (-1/2, 1), (0, 1/2), (0, 0), (1/2, 0), (1, -1/2), (1, 0), (3/2, 0), (1, 1/2). 
\end{displaymath}
By what we have already shown, there is a constant $C$ such that $\Xi(s, w) := ((s+10)(w+10))^{-C} (s-1)(w-1)(s+w-3/2)Z(s, w; \psi, \psi')$ is holomorphic and bounded in the tube $\{(s, w) \mid 4 < |\Re s|^2 + |\Re w|^2 < 5\}$.  A standard argument in several complex variables (see Propositions 4.6 and 4.7 and the argument on p.\ 341 of \cite{DGH}) shows that $\Xi(s, w)$ is holomorphic and bounded in the tube $\{(s, w) : |\Re s|^2 + |\Re w| < 5\}$. This completes the proof of the lemma.\\

Iterating \eqref{feqs1} and \eqref{feqs2} we find
\begin{equation}\label{feq3}
  \textbf{Z}(s, w) = \mathcal{B}(s) \cdot \mathcal{A} \cdot \mathcal{B}(s+w-1/2) \cdot \mathcal{A} \cdot \mathcal{B}(w) \cdot \mathcal{A} \, \textbf{Z}(1-s, 1-w)
\end{equation}

 
A computation shows that the matrix $\mathcal{M}(s, w) := \mathcal{B}(s) \cdot \mathcal{A} \cdot \mathcal{B}(s+w-1/2) \cdot \mathcal{A} \cdot \mathcal{B}(w) \cdot \mathcal{A}$ contains 124 zeros (out of 256 entries), but it is far from being diagonal. It would be nice to find a more symmetric version of \eqref{feq3}. 

An inspection of the matrix $\mathcal{B}$ in \eqref{defB} shows the following notationally more cumbersome, but slightly more practical form of \eqref{feq3}: There are absolute constants $\alpha^{(\text{\boldmath$\kappa$\unboldmath}, \textbf{j})}_{\rho, \rho', \psi, \psi' }$ such that
\begin{equation}\label{finalfunct}
\begin{split}
  Z(s, w; \psi, \psi')&  = \sum_{\kappa_1, \kappa_2, \kappa_3 = 0}^1 \sum_{j_1, j_2 = -6}^2 \sum_{\rho, \rho' \in \widehat{(\Bbb{Z}/8\Bbb{Z})^{\ast}}} \frac{\alpha^{(\text{\boldmath$\kappa$\unboldmath}, \textbf{j})}_{\rho, \rho', \psi, \psi' } 2^{j_1s + j_2w}\pi^{2s+2w-2}}{(4^s-4)(4^{s+w-1/2}-4)(4^w - 4)} \\
  & \times \frac{\Gamma(\frac{1-s+\kappa_1}{2})}{\Gamma(\frac{s+\kappa_1}{2})}\frac{\Gamma(\frac{3/2-s-w+\kappa_2}{2})}{\Gamma(\frac{s+w-1/2+\kappa_2}{2})}  \frac{\Gamma(\frac{1-w+\kappa_3}{2})}{\Gamma(\frac{w+\kappa_3}{2})}   Z(1-s, 1-w; \rho, \rho').
  \end{split}
\end{equation}

\section{A first approximate functional equation}

We use the functional equation \eqref{feq3} to obtain an explicit description of the function $Z(1/2+it, 1/2+iu; \psi, \psi')$, a so-called approximate functional equation. 
Although our assumptions are somewhat different, we follow essentially the     argument of \cite[Theorem 2.5]{Ha}. For $u, t \in \Bbb{R}$ we introduce henceforth  the following notation: let
 \begin{equation}\label{para}
  U := 1+|u|, \quad  T := 1+|t|, \quad S := 1 + |u+t|,  \quad X = STU
\end{equation}
and
\begin{displaymath}
   C  = 4 C(0, u) =  \left|\frac{1}{4}+\frac{i(u+t)}{2}\right| \cdot  \left|\frac{1}{4}+\frac{iu}{2}\right| 
 \end{displaymath}
with the notation  as in \eqref{anacon}. 

\begin{lemma}
There is a smooth, rapidly decaying function $V$,  and for any $u, t \in \Bbb{R}$   there are absolutely bounded  constants $\lambda^{\pm}_{j, \rho, \rho', \psi, \psi'}(u, t)$    such that for any $\varepsilon > 0$  and any $C' \geq C^{1/2+\varepsilon}$ one has\begin{displaymath}
\begin{split}
 & Z(1/2+it , 1/2 + iu; \psi, \psi') =  \sum_{\rho, \rho'} \sum_{j =-8}^{4} \sum_{\pm}\lambda^{\pm}_{j, \rho, \rho', \psi, \psi'}(u, t) \\
  & \times \sum_{\substack{d, m \text{ odd} \\dm^2 \leq C'}} \frac{L_2(1/2, \chi_d \rho) \rho'(d)}{(dm^2)^{1/2\pm iu}} V\left(\frac{dm^2}{2^j\sqrt{C}}\right) + O\left((TC)^{1/4+\varepsilon}\min(S, U)^{-1}\right).
  \end{split}
\end{displaymath}
 \end{lemma}

\textbf{Remark:} The error term can be improved with more careful estimations, but the above result suffices for our purposes. Note that $C'$ is bounded below, but otherwise independent of $u$ and $t$. \\

\textbf{Proof.} Let $t, u \in \Bbb{R}$.  Let $H$ be an even, holomorphic function with $H(0) = 1$ satisfying the growth estimate 
\begin{equation}\label{defH}
  H(z) \ll_{\Re z, A} (1+|z|)^{-A}
\end{equation}  
for any $A > 0$. Define 
\begin{displaymath}
  F_{u, t}(z) = \frac{1}{2}C^{-z/2} \frac{\Gamma(\frac{\frac{1}{2} - iu}{2})\Gamma(\frac{\frac{1}{2} - i(u+t)}{2})}{\Gamma(\frac{\frac{1}{2} + iu}{2})\Gamma(\frac{\frac{1}{2} + i(u+t)}{2})} \frac{\Gamma(\frac{\frac{1}{2} + iu+z}{2})\Gamma(\frac{\frac{1}{2} + i(u+t)+z}{2})}{\Gamma(\frac{\frac{1}{2} - iu-z}{2})\Gamma(\frac{\frac{1}{2} - i(u+t)-z}{2})}  +  \frac{1}{2}C^{z/2}. 
\end{displaymath}
Clearly $F_{u, t}$ is of moderate growth in fixed vertical strips and $F_{u, t}(0) = 1$. We consider the integral
\begin{equation}\label{afe1}
\begin{split}
  \frac{1}{2\pi i} \int_{(1)}& \frac{(2^{\frac{1}{2}+iu+z}-1)(2^{\frac{1}{2}+i(u+t)+z}-1)}{(2^{\frac{1}{2}+iu}-1)(2^{\frac{1}{2}+i(u+t)}-1)} \cdot 
  \frac{(4^{\frac{1}{2} + iu + z} - 4)(4^{\frac{1}{2} + i(u+t) + z} - 4)}{(4^{\frac{1}{2}+iu}-4)(4^{\frac{1}{2}+i(u+t)}-4)} \\
  & Z(1/2+it, 1/2 + iu + z; \psi, \psi') F_{u, t}(z) H(z) \frac{dz}{z}.
  \end{split}
\end{equation}
The first fraction cancels the possible poles at $z=-1/2 - iu$, $z = -1/2-iu-it$ of $F_{u, t}$, the second fraction cancels the possible poles at $z=1/2 - iu$, $z = 1/2-iu-it$ of $Z$ and goes well with the functional equation \eqref{finalfunct}. This device is not strictly necessary, but it is convenient. 
We shift the contour to $\Re z = -1$. The pole at $z=0$ contributes 
\begin{equation}\label{afe2}
    Z(1/2+it, 1/2 + iu; \psi, \psi').  
\end{equation}
In the remaining integral we apply the functional equation \eqref{finalfunct} together with \eqref{cot} and make a change of variables $z \mapsto -z$ getting
\begin{equation}\label{afe3}
\begin{split}
&  - \frac{1}{2\pi i} \int_{(1)}  \frac{(2^{\frac{1}{2}+iu-z}-1)(2^{\frac{1}{2}+i(u+t)-z}-1)}{(2^{\frac{1}{2}+iu}-1)(2^{\frac{1}{2}+i(u+t)}-1)}   
  \sum_{\kappa_1, \kappa_2, \kappa_3 = 0}^1 \sum_{j_1, j_2 = -6}^2 \sum_{\rho, \rho' \in \widehat{(\Bbb{Z}/8\Bbb{Z})^{\ast}}}  \alpha^{(\text{\boldmath$\kappa$\unboldmath}, \textbf{j})}_{ \rho, \rho', \psi, \psi' } \\
  & \times \frac{2^{j_1(\frac{1}{2}+it) + j_2(\frac{1}{2} + iu - z)}\pi^{2i(u+t)-2z}}{(4^{\frac{1}{2}+it}-4)(4^{\frac{1}{2}+iu}-4)(4^{\frac{1}{2}+i(u+t)}-4)}  \frac{\Gamma(\frac{\frac{1}{2} - it + \kappa_1}{2})\Gamma(\frac{\frac{1}{2} - i(u+t) + z}{2})\Gamma(\frac{\frac{1}{2} - iu + z}{2})}{\Gamma(\frac{\frac{1}{2} + it +\kappa_1 }{2})\Gamma(\frac{\frac{1}{2} + i(u+t) -z  }{2})\Gamma(\frac{\frac{1}{2} + iu -z  }{2})} \\
  &  \times \cot\left(\frac{\pi(\frac{1}{2} + i(u+t) - z)}{2}\right)^{\kappa_2}   \cot\left(\frac{\pi(\frac{1}{2} + iu - z)}{2}\right)^{\kappa_3}  \\
  & \times  Z(1/2-it , 1/2 - iu + z; \rho, \rho') F_{u, t}(-z) H(z) \frac{dz}{z}.
\end{split}
\end{equation}
Then \eqref{afe1} equals the sum of \eqref{afe2} and \eqref{afe3}. We need to simplify the unduly complicated term \eqref{afe3}. We observe that
\begin{displaymath}
   \frac{\Gamma(\frac{\frac{1}{2} - iu + z}{2})\Gamma(\frac{\frac{1}{2} - i(u+t) + z}{2})}{\Gamma(\frac{\frac{1}{2} + iu -z  }{2})\Gamma(\frac{\frac{1}{2} + i(u+t) -z  }{2})} F_{u, t}(-z)   = \frac{\Gamma(\frac{\frac{1}{2} - iu}{2})\Gamma(\frac{\frac{1}{2} - i(u+t)}{2})}{\Gamma(\frac{\frac{1}{2}+iu}{2})\Gamma(\frac{\frac{1}{2}+i(u+t)}{2})} F_{-u, -t}(z), 
\end{displaymath}
so that \eqref{afe3} simplifies to
\begin{equation}\label{afe4}
  \begin{split}
 & - \frac{1}{2\pi i} \int_{(1)}  \sum_{\kappa_2, \kappa_3 = 0}^1 \sum_{j = -6}^4 \sum_{\rho, \rho' \in \widehat{(\Bbb{Z}/8\Bbb{Z})^{\ast}}} \mu^{(\kappa_2, \kappa_3, j)}_{ \rho, \rho', \psi, \psi' }(u, t) 2^{-jz} \pi^{-2z}    \cot\left(\frac{\pi(\frac{1}{2} + i(u+t) - z)}{2}\right)^{\kappa_2}\\
 & \times   \cot\left(\frac{\pi(\frac{1}{2} + iu - z)}{2}\right)^{\kappa_3}  Z(1/2-it, 1/2 - iu + z; \rho, \rho') F_{-u, -t}(z) H(z) \frac{dz}{z}
\end{split}
\end{equation}
for certain absolutely bounded complex numbers $\mu^{(\kappa_2, \kappa_3, j)}_{ \rho, \rho', \psi, \psi' }(u, t)$. 
In \eqref{afe1} and \eqref{afe4} we open the Dirichlet series using the definition \eqref{defZ}. This yields   the following preliminary version of the lemma: for $\kappa_2, \kappa_3 \in \{0, 1\}$ and $u, t \in \Bbb{R}$ let
\begin{equation}\label{defV}
\begin{split}
  V_{u, t}^{( \kappa_2, \kappa_3)}(\xi) = \frac{1}{2 \pi i} \int_{(1)}& \cot\left(\frac{\pi(\frac{1}{2} + i(u+t) - z)}{2}\right)^{\kappa_2} \cot\left(\frac{\pi(\frac{1}{2} + iu - z)}{2}\right)^{\kappa_3}  \\\
  & \times C^{-z/2} F_{ u, t}(z)\pi^{-2z}  H(z) \xi^{-z} \frac{dz}{z}.  
\end{split}
\end{equation}
Then there are absolutely bounded constants $\mu^{\pm, (\kappa_2, \kappa_3, j)}_{ \rho, \rho', \psi, \psi'}(u, t) \in \Bbb{C}$ such that
\begin{equation}\label{linearcombi}
\begin{split}
  Z(1/2+it, 1/2+iu; \psi, \psi') & =  \sum_{\rho, \rho'} \sum_{j =-6}^4 \sum_{\kappa_2, \kappa_3 = 0}^1 \sum_{\pm} \mu^{\pm, (\kappa_2, \kappa_3, j)}_{ \rho, \rho', \psi, \psi'}(u, t)\\
  & \sum_{d, m\text{ odd}} \frac{L_2(1/2+it, \chi_d \rho) \rho'(d)}{d^{\frac{1}{2} \pm iu}m^{1\pm 2i(u+t)}} V^{(\kappa_2, \kappa_3)}_{\pm u, \pm t }\left(\frac{2^jdm^2}{\sqrt{C}}\right).  
    \end{split}
\end{equation}
We analyze the function $V_{u, t}^{(\kappa_2, \kappa_3)}$ and  quote two bounds of \cite{Ha} (see also the erratum): By Lemma 3.1 in \cite{Ha} we have
\begin{equation}\label{Ha1}
  C^{-z/2} F_{u, t}(z)  \ll (1+|z|)^{2 \Re z}, \quad \Re z \geq 0
\end{equation}
and by Lemma 4.1 in \cite{Ha} we have
\begin{equation}\label{Ha2}
  C^{-\frac{iy}{2}} F_{u, t}(iy) - 1 \ll |y| C^{\varepsilon}\min(S, U)^{-1}, \quad y \in \Bbb{R},  |y| < C^{\varepsilon}
\end{equation}
for any $0 < \varepsilon < 1/2$. Both \eqref{Ha1} and \eqref{Ha2} are uniform in $u$ and $t$. 

Now we return to \eqref{defV} and shift the contour to the far right. We pick up possible poles at $z = 5/2+2n  -iu$ and $z=5/2+2n-i(u+t)$, $n=0, 1, 2, \ldots$ from the cotangent, whose contribution is by \eqref{defH} and \eqref{Ha1} at most $\ll \log(2+\xi) \xi^{-5/2} \min(S, U)^{-A}$ (we need the logarithm if $t=0$). Hence we find
\begin{displaymath}
  V^{(\kappa_2, \kappa_3)}_{u, t}(\xi) \ll  \log(2+\xi) \xi^{-5/2} \min(S, U)^{-A} +\xi^{-A},
\end{displaymath}  
uniformly in $u, t$. Combining this with the average bound \eqref{HBest}, we conclude
\begin{equation}\label{linearcombi1}
\begin{split}
 & Z(1/2+it, 1/2+iu; \psi, \psi')  =  \sum_{\rho, \rho'} \sum_{j =-6}^4 \sum_{\kappa_2, \kappa_3 = 0}^1 \sum_{\pm} \mu^{\pm, (\kappa_2, \kappa_3, j)}_{\rho, \rho', \psi, \psi'}(u, t)\\
  & \times \sum_{\substack{d, m\text{ odd}\\ dm^2 \leq C^{1/2+\varepsilon}}} \frac{L_2(1/2+it, \chi_d \rho) \rho'(d)}{d^{\frac{1}{2} \pm iu}m^{1\pm 2i(u+t)}} V^{(\kappa_2, \kappa_3)}_{\pm u, \pm t }\left(\frac{2^jdm^2}{ \sqrt{C}}\right) +O\left((TC)^{1/4+\varepsilon}\min(S, U)^{-A}\right)
  \end{split}
\end{equation}
for any $\varepsilon > 0$.  We can now remove the dependence on $u$ and $t$ of $V_{u, t}^{(\kappa_2, \kappa_3)}$. To this end we consider
\begin{equation}\label{newV}
   V (\xi) = \frac{1}{2 \pi i} \int_{(1)}  \pi^{-2z} H(z) \xi^{-z} \frac{dz}{z}.  
\end{equation}
This a smooth, rapidly decreasing function.  For $\xi \ll C^{\varepsilon}$ we estimate the difference 
\begin{equation}\label{diff}
\begin{split}
 & V_{u, t}^{(\kappa_2, \kappa_3)}(\xi) - (-i \, \text{sign} (u+t))^{\kappa_2}(-i \, \text{sign} (u))^{\kappa_3}V(\xi)\\
  & = \frac{1}{2 \pi i} \int_{\gamma}  \left(\cot\left(\frac{\pi(\frac{1}{2} + i(u+t) - z)}{2}\right)^{\kappa_2}  \cot\left(\frac{\pi(\frac{1}{2} + iu - z)}{2}\right)^{\kappa_3} C^{-z/2} F_{u, t}(z)\right.\\
  &\quad\quad\Bigl. - (-i \, \text{sign} (u+t))^{\kappa_2}(-i \, \text{sign} (u))^{\kappa_3} \Bigr)  
  \pi^{-2z} H(z) \xi^{-z} \frac{dz}{z}
  \end{split}
\end{equation}
where $\gamma = \gamma_1\gamma_2\gamma_3$ with $\gamma_1 = [-i\infty, -i\varepsilon]$, $\gamma_2$   a semicircle to the right of the origin joining $-i\varepsilon $ and $i\varepsilon$, and $\gamma_3 = [i\varepsilon, i\infty]$. The portion $|\Im z| > C^{\varepsilon}$ contributes by \eqref{Ha1} and \eqref{defH} at most $O(C^{-A})$. In the remaining part we insert the formula \eqref{cotbound} at the cost of an error $O(C^{\varepsilon} \min(S, U)^{-A})$. The integrand is now holomorphic at $z=0$, and we replace the semicircle $\gamma_2$ with a straight line through the origin. Now we insert \eqref{Ha2} and bound the integral \eqref{diff} by $O(C^{\varepsilon}\min(S, U)^{-1}) $. 

Now we  replace  in \eqref{linearcombi1} the weight function $V_{u, t}^{(\kappa_2, \kappa_3)}$ by $(-i \, \text{sign} (u+t))^{\kappa_2}(-i \, \text{sign} (u))^{\kappa_3}V $ at the cost of an error
\begin{displaymath}
  C^{\varepsilon}\min(S, U)^{-1} \sum_{dm^2 \leq C^{\frac{1}{2} + \varepsilon}} \frac{|L_2(1/2+it, \chi_d \rho)|}{d^{1/2}m}  \ll  (TC)^{1/4+\varepsilon}\min(S, U)^{-1}
\end{displaymath}
by \eqref{HBest}. The lemma follows now with $C' = C^{1/2 +\varepsilon}$. By the rapid decay of $V$ it remains valid for any larger $C'$.  \\

Lemma 3 reduces the estimation of $Z(1/2+it, 1/2+iu; \psi, \psi')$ to bounding
\begin{displaymath}
\sum_{\substack{d, m \text{ odd} \\dm^2 \leq C'}} \frac{L_2(1/2+it, \chi_d \rho) \rho'(d)}{d^{1/2\pm iu}m^{1\pm 2i(u+t)}} V\left(\frac{2^jdm^2}{\sqrt{C}}\right).
\end{displaymath}
Applying a smooth partition of unity, it is therefore enough to bound
\begin{displaymath}
  D_{\psi, \psi'}(t, u, P; W) := \sum_{d, m \text{ odd}} \frac{L_2(1/2 + it, \chi_d\psi) \psi'(d)}{d^{1/2 + iu}m^{1+2i(u+t)}} W\left(\frac{dm^2}{P}\right) 
\end{displaymath}
for a smooth function $W$ with support on $[1, 2]$ and 
\begin{equation}\label{defP}
  1\leq P \leq (US)^{1/2+\varepsilon}. 
\end{equation}
Henceforth we will always assume that $P$ satisfies \eqref{defP}, and we recall the notation \eqref{para}. We prove the following variant of the preceding lemma.
\begin{lemma} Let $\delta_0$ be given by \eqref{fepara} and let $T' \geq TX^{\varepsilon}$. Then the  following bound holds:
\begin{displaymath}
\begin{split}
  D_{\psi, \psi'}(t, u, P; W) \ll  &  \sum_{\pm} \sum_{m \leq P^{1/2+\varepsilon}} \frac{(Xm)^{\varepsilon}}{m}  \int_{\varepsilon-iX^{\varepsilon}}^{\varepsilon+iX^{\varepsilon}} \int_{\varepsilon-iX^{\varepsilon}}^{\varepsilon+iX^{\varepsilon}} \\
  &  \Biggl|  \sum_{\substack{d_0 \leq m^{-2}P^{1+\varepsilon}\\  d_0 \text{ odd}}}\sum_{n \leq (T'P)^{1/2}} \frac{(\chi_{d_0}\psi)(n) \psi'(d_0) \delta_0^{s/2}}{n^{1/2 \pm it-s}d_0^{1/2 + iu-w}}    \Biggr| \, |dw|\, |ds| +P^{-A}.
  \end{split}
\end{displaymath}
\end{lemma}


\textbf{Proof.} This follows quickly from \eqref{approxL} and \eqref{HBbilinear}. More precisely,   by \eqref{nonfund} and \eqref{approxL} we have
\begin{displaymath}
 | D_{\psi, \psi'}(t, u, P; W)| \ll \sum_{\pm} \sum_{m, d_1}\frac{d_1^{\varepsilon}}{md_1}\Biggl| \sum_{  d_0 \text{ odd}}\sum_{n} \frac{(\chi_{d_0}\psi)(n) \psi'(d_0)}{n^{1/2 \pm it}d_0^{1/2 + iu}} G^{(\psi)}_{\pm t}\left(\frac{n}{\sqrt{\delta_0}}\right) W\left(\frac{d_0d_1^2m^2}{P}\right)\Biggr|
\end{displaymath}
where 
$G^{(\psi)}_{t}$ is given by \eqref{defG}. By the rapid decay of $W$ and $G_t$ (cf.\ \eqref{defG}) we can truncate the sums at  $d_0d_1^2m^2 \leq P^{1+\varepsilon}$ and $n \leq (T'P)^{1/2 }$ at the cost of an error $O(P^{-A})$.  Let $\widehat{W}$ denote the Mellin transform of $W$. Then $\widehat{W}$ is an entire function with rapid decay in fixed vertical strips. We recast $W$ and $G_{t}^{(\psi)}$ by Mellin inversion getting
\begin{displaymath}
\begin{split}
 | D_{\psi, \psi'}(t, u, P; W)| \ll & P^{-A} + \sum_{\pm} \sum_{m \leq P^{1/2+\varepsilon}} \frac{1}{m^{1-\varepsilon}}  \int_{(\varepsilon)} \int_{(\varepsilon)} \Biggl| \left(\cos\frac{\pi s}{4A}\right)^{-4A} \frac{\Gamma(\frac{\frac{1}{2} + it +s +\kappa}{2})}{\Gamma(\frac{\frac{1}{2} + it +\kappa}{2})}\Biggr. \\
 & \times  \Biggl. \sum_{\substack{d_0 \leq m^{-2}P^{1+\varepsilon}\\  d_0 \text{ odd}}}\sum_{n \leq (TP)^{1/2+\varepsilon}} \frac{(\chi_{d_0}\psi)(n) \psi'(d_0) \delta_0^{s/2}}{n^{1/2 \pm it-s}d_0^{1/2 + iu-w}}   \widehat{W}(w) \left(\frac{P}{m^2}\right)^w\Biggr|\, |dw|\, \left|\frac{ds}{s}\right|
 \end{split}
\end{displaymath}
with $\kappa$ as in \eqref{fepara}. By Stirling's formula \eqref{stir2} and the rapid decay of   $\widehat{W}$ we can truncate the $s, w$-integration, and the lemma follows. 

\section{A second approximate functional equation}

In this section we establish a different bound for $D_{\psi, \psi'}(t, u, P;W)$. By Mellin inversion we have 
\begin{equation}\label{hence}
  D_{\psi, \psi'}(t, u, P; W) = \frac{1}{2\pi i} \int_{(1)} Z(1/2+it , 1/2 + iu + w; \psi, \psi') \widehat{W}(w) P^w dw.   
\end{equation}
Let $H$ be the same function as in the preceding proof, that is,  $H$ is even and holomorphic, rapidly decaying in fixed vertical strips and $H(0) = 1$. For $\Re z \geq 3/2$ and $R>0$  we consider the term
\begin{equation}\label{expression}
- Z(1/2+it, z; \psi, \psi') + \frac{1}{2\pi i} \int_{(3)} \frac{4^{\frac{1}{2}+it+ s}-4}{4^{\frac{1}{2}+it}-4} Z(1/2 +it+ s, z; \psi, \psi') R^s H(s) \frac{ds}{s}.
\end{equation}
We shift the contour to $\Re s = -3$. The possible pole of $Z$ at $s=1/2-it$ is cancelled by the first fraction, and the possible pole at $s=1-it- z$  contributes at most  $O(R^{1-\Re z}(1+|z+it|)^{-A})$. We change variables $s \mapsto -s$ and apply (one component of) the functional equation \eqref{feqs2}. Hence the preceding expression equals
\begin{equation}\label{expression2}
\begin{split}
  -\frac{1}{2\pi i}  & \int_{(3)}  \sum_{\kappa  = 0}^1 \sum_{j=-3}^1 \sum_{ \rho, \rho' \in \widehat{(\Bbb{Z}/8\Bbb{Z})^{\ast}}} \tilde{\alpha}^{(\kappa, j)}_{ \rho, \rho', \psi, \psi'}(t) 2^{-js} \pi^{-s}  \cot\left(\frac{\pi(\frac{1}{2}-it -s)}{2}\right)^{\kappa } \\
  & \times \frac{\Gamma(\frac{\frac{1}{2}+it+s}{2})}{\Gamma(\frac{\frac{1}{2}-it-s}{2})} R^{-s/2} Z(1/2 -it+ s, z+it-s; \rho, \rho') H(s) \frac{ds}{s} + O(R^{1-\Re z}(1+|it+z|)^{-A})
  \end{split}
\end{equation}
for absolutely bounded constants $\tilde{\alpha}^{\kappa, j}_{\rho, \rho', \psi, \psi'}(t)$. We substitute \eqref{expression} and \eqref{expression2} with $R = (PT)^{1/2}$   and $z = 1/2 + iu + w$ into \eqref{hence} getting
\begin{equation}\label{mainquant}
  D_{\psi, \psi'}(t, u, P; W) = D + \tilde{D} + O(P^{3/4}S^{-A}),
\end{equation}
say, where
\begin{equation}\label{defD}
\begin{split}
  D :=  \left(\frac{1}{2 \pi i}\right)^2 \int_{(1)}&\int_{(3)}  \frac{4^{\frac{1}{2}+it+s}-4}{4^{\frac{1}{2}+it}-4} Z(1/2+it+s, 1/2+iu+w; \psi, \psi')\widehat{W}(w)\frac{H(s)}{s} T^{\frac{s}{2}}P^{w+\frac{s}{2}} ds\, dw
  \end{split}
\end{equation}
and
\begin{equation}\label{defDD}
\begin{split}
 \tilde{D} & :=  \sum_{\kappa  = 0}^1 \sum_{j=-3}^1 \sum_{ \rho, \rho' \in \widehat{(\Bbb{Z}/8\Bbb{Z})^{\ast}}} \tilde{\alpha}^{(\kappa, j)}_{ \rho, \rho', \psi, \psi'}(t)  \left(\frac{1}{2 \pi i}\right)^2 \int_{(1)}\int_{(3)} \frac{\Gamma(\frac{\frac{1}{2}+it+s}{2})}{\Gamma(\frac{\frac{1}{2}-it-s}{2})}   \cot\left(\frac{\pi(\frac{1}{2}-it -s)}{2}\right)^{\kappa } \\
& \times  2^{-js} \pi^{-s} Z(1/2-it +s, 1/2+i(u+t)+w-s; \rho, \rho')\widehat{W}(w)   \frac{H(s)}{s} T^{-\frac{s}{2}} P^{w-\frac{s}{2}} ds\, dw. 
 \end{split}
\end{equation}
For the error term in \eqref{mainquant} we used the rapid decay of $\widehat{W}$. 
Both double integrals are absolutely convergent. In \eqref{defD}, we shift the $w$-integration to $\Re w = -1$ and change variables $w \mapsto -w$. There is a  possible pole on the way at $w=1/2-iu$   whose contribution is, by the rapid decay of $\widehat{W}$, bounded by $O(U^{-A})$. Hence
\begin{equation}\label{defDa}
\begin{split}
D =  \left(\frac{1}{2 \pi i}\right)^2 \int_{(1)}&\int_{(3)}  \frac{4^{\frac{1}{2}+it+s}-4}{4^{\frac{1}{2}+it}-4} Z(1/2+it+s, 1/2+iu-w; \psi, \psi')\\
  & \times\widehat{W}(-w)\frac{H(s)}{s}T^{\frac{s}{2}}P^{-w+\frac{s}{2}} ds\, dw + O(U^{-A}). 
\end{split}
\end{equation}
Here it is important to note that the partition of unity has removed the pole at $w=0$ that  would occur  if $\widehat{W}(w)$ was replaced by $H(w)/w$ as in Lemma 3. 

Now we apply  the functional equations \eqref{feqs1} and \eqref{feqs2} in the form
\begin{equation}\label{newfe1}
  \textbf{Z}(1/2 + it+ s, 1/2 + iu - w) =  \mathcal{A} \cdot \mathcal{B}(1/2 + iu - w)    \,\textbf{Z}(1/2 - iu+w, 1/2+i(u+t)-w+s)
\end{equation}
and 
\begin{equation}\label{newfe2}
  \textbf{Z}(1/2 -it + s, 1/2 + i(u+t) + w - s) = \mathcal{A} \cdot \mathcal{B}(1/2 + i(u+t) + w -s) \cdot \mathcal{A} \,\textbf{Z}(1/2 + iu + w, 1/2 - i(u+t) - w + s).
\end{equation}
For convenience we write this out explicitly: The $(\psi, \psi')$-component of \eqref{newfe1} is 
\begin{equation}\label{newfe1a}
\begin{split}
&  Z(1/2 + it+ s, 1/2 + iu - w; \psi, \psi') = \sum_{\kappa = 0}^1 \sum_{j = -3}^1 \sum_{\rho, \rho' \in \widehat{(\Bbb{Z}/8\Bbb{Z})^{\ast}}}  \beta^{(\kappa,  j)}_{\rho, \rho', \psi, \psi' }(u)\frac{2^{-j w}\pi^{-w}}{4^{\frac{1}{2}+iu-w} - 4} \\
&\times\frac{\Gamma(\frac{\frac{1}{2}-iu+w}{2})}{\Gamma(\frac{\frac{1}{2}+iu-w}{2})} \cot\left(\frac{\pi(\frac{1}{2} +iu-w)}{2}\right)^{\kappa}   Z(1/2-iu+w, 1/2+i(u+t)-w+s; \rho, \rho')
\end{split}
\end{equation}
for certain absolutely bounded constants $\beta^{(\kappa, j)}_{\rho, \rho', \psi, \psi'}(u)$, and the $(\rho, \rho')$-component of \eqref{newfe2} is 
\begin{equation}\label{newfe2a}
\begin{split}
&  Z(1/2 -it + s, 1/2 + i(u+t) + w-s; \rho, \rho') = \sum_{\kappa = 0}^1 \sum_{j = -3}^1 \sum_{\tilde{\rho}, \tilde{\rho}' \in \widehat{(\Bbb{Z}/8\Bbb{Z})^{\ast}}}  \tilde{\beta}^{(\kappa,  j)}_{ \tilde{\rho}, \tilde{\rho}', \rho, \rho' }(u, t)\frac{2^{j(w-s)}\pi^{w-s}}{4^{\frac{1}{2}+i(u+t)+w-s} - 4} \\
&\times\frac{\Gamma(\frac{\frac{1}{2}-i(u+t)-w+s}{2})}{\Gamma(\frac{\frac{1}{2}+i(u+t)+w-s}{2})} \cot\left(\frac{\pi(\frac{1}{2} +i(u+t)+w-s)}{2}\right)^{\kappa}   Z(1/2+iu+w, 1/2-i(u+t)-w+s; \tilde{\rho}, \tilde{\rho}')
\end{split}
\end{equation} 
for certain absolutely bounded constants $\tilde{\beta}^{(\kappa, j)}_{ \tilde{\rho}, \tilde{\rho}', \rho, \rho'}(u, t)$. We substitute \eqref{newfe1a} into \eqref{defDa} and open both components of the absolutely convergent double Dirichlet series. In this way we see
\begin{equation}\label{finalD}
\begin{split}
  D & = \sum_{\rho, \rho' \in \widehat{(\Bbb{Z}/8\Bbb{Z})^{\ast}}}  \sum_{n, d \text{ odd}} \frac{\chi_d(n)\rho(n)\rho'(d)}{n^{1/2-iu} d^{1/2+i(u+t)}} \,\,   \left(\frac{1}{2 \pi i}\right)^2 \int_{(1)}\int_{(3)}\zeta_2(2s+2it+ 1) \\
  & \times \sum_{\kappa = 0}^1 \sum_{j = -3}^1 \beta^{(\kappa,  j)}_{ \rho, \rho', \psi, \psi' }(u)  \frac{4^{\frac{1}{2}+it+s}-4}{4^{\frac{1}{2}+it}-4}   \frac{2^{-jw}\pi^{-w}}{4^{\frac{1}{2}+iu-w} - 4} \cot\left(\frac{\pi(\frac{1}{2} +iu-w)}{2}\right)^{\kappa} \frac{\Gamma(\frac{\frac{1}{2}-iu+w}{2})}{\Gamma(\frac{\frac{1}{2}+iu-w}{2})}  \\
&  \times \left(\frac{nP}{d}\right)^{-w} \Biggl(\frac{d}{\sqrt{T P}}\Biggr)^{-s} \widehat{W}(-w) \frac{H(s)}{s}  ds\, dw + O(U^{-A}).
  \end{split} 
\end{equation}
Similarly, we substitute \eqref{newfe2a} into \eqref{defDD} getting
\begin{equation}\label{finalDD}
\begin{split}
 &\tilde{D}  = \sum_{\rho, \rho' \in \widehat{(\Bbb{Z}/8\Bbb{Z})^{\ast}}}  \sum_{n, d \text{ odd}} \frac{\chi_d(n)\rho(n)\rho'(d)}{n^{1/2+iu} d^{1/2-i(u+t)}} \,\,   \left(\frac{1}{2 \pi i}\right)^2 \int_{(1)}\int_{(3)} \sum_{\kappa_1, \kappa_2 = 0}^1 \sum_{j_1, j_2 = -3}^1 \gamma^{(\kappa_1,  \kappa_2, j_1, j_2)}_{ \rho, \rho', \psi, \psi' } (u, t) \\
 &\zeta_2(2s+2it+ 1)\frac{2^{-j_1s + j_2(  w-s)}\pi^{w-2s}}{4^{\frac{1}{2}+i(u+t)+w-s} - 4}\cot\left(\frac{\pi(\frac{1}{2}-it -s)}{2}\right)^{\kappa_1 }  \cot\left(\frac{\pi(\frac{1}{2} +i(u+t)+w-s)}{2}\right)^{\kappa_2} \\
 & \frac{\Gamma(\frac{\frac{1}{2}+it+s}{2})}{\Gamma(\frac{\frac{1}{2}-it-s}{2})}   \frac{\Gamma(\frac{\frac{1}{2}-i(u+t)-w+s}{2})}{\Gamma(\frac{\frac{1}{2}+i(u+t)+w-s}{2})}  \widehat{W}(w)       \frac{H(s)}{s} \left(\frac{n}{dP}\right)^{-w} \left(d\sqrt{T P}\right)^{-s} ds\, dw
\end{split}
\end{equation}
for certain absolutely bounded constants $ \gamma^{(\kappa_1,  \kappa_2, j_1, j_2)}_{ \rho, \rho', \psi, \psi' }(u, t)$. We can now substitute \eqref{finalD} and \eqref{finalDD} into \eqref{mainquant} to obtain an approximate functional equation for the quantity $D_{\psi, \psi'}(t, u, P; W)$ that we want to bound. It remains to analyze the two double integrals in \eqref{finalD} and \eqref{finalDD} as functions of $d, n, t, u, P$. To this end we shift contours and use Stirling's formula \eqref{stir}. Let 
\begin{equation}\label{newpara}
  T' \geq TX^{\varepsilon}, \quad S' \geq SX^{\varepsilon}, \quad U' \geq UX^{\varepsilon}, \quad X':= S'T'U'.
\end{equation}  
   First we shift in \eqref{finalD} the $s$-contour to $\Re s = A$, bounding the double integral by $\ll \frac{dU}{nP} \left(\frac{\sqrt{TP}}{d}\right)^A$. Hence we can restrict the $d$-summation to
\begin{displaymath}
  d \leq (T' P)^{1/2} 
\end{displaymath}
with an error of at most $O(X^{-A})$. Next we shift the $w$-integration to $\Re w = A$. The poles of the cotangent are cancelled by the poles of the Gamma function. On $\Re s = \Re w = A$ we bound the double integral by $\ll (UT)/(\sqrt{P} n)^A$, hence we can  we can truncate the $n$-sum at 
\begin{displaymath}
  n \leq  \frac{(T')^{1/2}U'}{P^{1/2}}
\end{displaymath}
with the same  error. Similarly,  in \eqref{finalDD} we shift the $s$-contour to $\Re s = A$ and restrict the $d$-sum to 
\begin{displaymath}
  d \leq \frac{(T')^{1/2}S'}{P^{1/2}},
\end{displaymath}  
at the cost of an error $O(X^{-A})$; then we shift the $w$-contour to $\Re w = A-1$ and truncate the $n$-sum at 
\begin{displaymath}   
  n \leq   (T'P)^{1/2}, 
\end{displaymath}
again with an error $O(X^{-A})$. Having truncated the double sums \eqref{finalD} and \eqref{finalDD} in this way, we shift the contours back to $\Re s = \Re w = \varepsilon$ and interchange the (now finite) $d, n$-double sum with the absolutely convergent $s, w$-double integral. Finally, by the rapid decay of $\widehat{W}$ and $H$ we can truncate the $s, w$-integration. Thus we arrive at the following complement to Lemma 4. We keep the notation \eqref{para} and \eqref{newpara}. 
\begin{lemma}
The following bound holds:
\begin{displaymath}
\begin{split}
  D_{\psi, \psi'}(t, u, P; W)  \ll (X')^{\varepsilon} & \max_{\rho, \rho' \in \widehat{(\Bbb{Z}/8\Bbb{Z})^{\ast}}}  \int_{\varepsilon-i(X')^{\varepsilon}}^{\varepsilon+i(X')^{\varepsilon}} \int_{\varepsilon-i(X')^{\varepsilon}}^{\varepsilon+i(X')^{\varepsilon}}   \Biggl|\sum_{\substack{n, d \text{ odd}\\ d\leq (T'P)^{1/2}  \\ n \leq (T'/P)^{1/2}U' }} \frac{\chi_d(n)\rho(n)\rho'(d)}{n^{1/2-iu+w} d^{1/2+i(u+t)-w+s}}\Biggr|  \\
  & +     \Biggl| \sum_{\substack{n, d \text{ odd}\\ d\leq (T'/P)^{1/2}S'  \\ n \leq (T'P)^{1/2}  }} \frac{\chi_d(n)\rho(n)\rho'(d)}{n^{1/2+iu+w} d^{1/2-i(u+t)-w+s}}\Biggr| \, |ds\, dw|   +P^{3/4}S^{-A} .
   \end{split}
\end{displaymath}
\end{lemma}



\section{Proof of Theorem 1 and 2}

We are now prepared to prove our main theorems.    In order to prove Theorem 1, we first observe that without loss of generality we can assume 
\begin{equation}\label{wlog}
  T \leq U \asymp S.
\end{equation}
Indeed, by \eqref{feqs1} and \eqref{feqs2} we have
\begin{displaymath}
\begin{split}
 &  Z(1/2 + it, 1/2+iu; \psi, \psi') \ll \max_{\rho, \rho'} |Z(1/2+iu, 1/2+it; \rho, \rho')|,\\
 &  Z(1/2 + it, 1/2+iu; \psi, \psi') \ll \max_{\rho, \rho'} |Z(1/2+i(t+u), 1/2-iu; \rho, \rho')|
\end{split}
\end{displaymath}
with absolute implied constants. Hence we can exchange $u$ and $t$, if necessary, to ensure $|t| \leq |u|$, and then we can exchange $t$ and $t+u$ (thereby sending $u$ to $-u$), if necessary, to ensure $|u|/2 \leq |u+t| \leq 2|u|$. The desired bound of Theorem 1 is symmetric in these permutations which justifies our assumption \eqref{wlog}. 

 In Lemma 4 and 5 we estimate the character sum by \eqref{HBbilinear} and conclude that
\begin{displaymath}
  D_{\psi, \psi'}(t, u, P; W) \ll U^{\varepsilon} \min \left(P^{1/2} + (TP)^{1/4}, (TP)^{1/4} +\left(\frac{T}{P}\right)^{1/4} U^{1/2}\right).
\end{displaymath}
Lemma 3 implies now
\begin{displaymath}
  Z(1/2+it, 1/2+iu, \psi, \psi') \ll \frac{(TU^2)^{1/4+\varepsilon}}{U} + U^{\varepsilon} \max_{P \ll U}\left( (TP)^{1/4} + \min \left(P^{1/2} ,  \left(\frac{T}{P}\right)^{1/4} U^{1/2}\right)\right).
\end{displaymath}
If $P \leq U^{2/3}T^{1/3}$ we take the first term in the parentheses, otherwise the second. In either case, 
\begin{displaymath}
   Z(1/2+it, 1/2+iu, \psi, \psi') \ll U^{2/3+\varepsilon}T^{1/3}
\end{displaymath}
what was to be proved. \\

We proceed to prove Theorem 2. Let $W$ be a nonnegative function satisfying $W(x) = 1$ for $x \in [-1, 1]$ and $W(x) = 0$ for $|x| \geq 2$. We need to prove
\begin{equation}\label{toprove}
   \int \int W\left(\frac{t}{Y_1}\right) W\left(\frac{u}{Y_2}\right)  |Z(1/2+it, 1/2+iu; \psi, \psi')|^2 du \, dt \ll   (Y_1Y_2)^{1+\varepsilon}.
\end{equation}
 The same argument as above shows that without loss of generality we can assume 
 \begin{equation}\label{assump}
   Y_1 \leq Y_2 
 \end{equation}
     Let $\gamma$ be the vertical segment  $[\varepsilon - iY_2^{\varepsilon},  \varepsilon + iY_2^{\varepsilon}]$. For $P \leq Y_2^{1+\varepsilon}$ let 
  \begin{displaymath}
\begin{split}
& Q^{(P, Y_1)}_{1, \pm}(t, u; s, w) :=    \sum_{\substack{D = 2^{\nu_1} \leq P^{1+\varepsilon}\\ N = 2^{\nu_2} \leq (Y_1P)^{1/2+\varepsilon}}} \Biggl|  \sum_{\substack{ d_0 \text{ odd}\\ d_0 \sim D, \, n\sim N }}  \frac{(\chi_{d_0}\psi)(n) \psi'(d_0) \delta_0^{s/2}}{n^{1/2 \pm it-s}d_0^{1/2 + iu-w}}    \Biggr| ,\\
& Q^{(P, Y_1, Y_2)}_2(t, u; s, w) :=  \sum_{\substack{D = 2^{\nu_1} \leq (Y_1P)^{1/2+\varepsilon}\\N = 2^{\nu_2} \leq (\frac{Y_1}{P})^{1/2}Y_2^{1+\varepsilon}}}   \Biggl| \sum_{\substack{n, d \text{ odd}\\ d\sim D, n \sim N}} \frac{\chi_d(n)\rho(n)\rho'(d)}{n^{1/2- iu+w} d^{1/2+ i(u+t)-w+s}}\Biggr|,\\
&  Q^{(P, Y_1, Y_2)}_3(t, u; s, w) :=  \sum_{\substack{D = 2^{\nu_1} \leq  (\frac{Y_1}{P})^{1/2}Y_2^{1+\varepsilon} \\ N=2^{\nu_2} \leq (Y_1P)^{1/2} Y_2^{\varepsilon}}}    \Biggl|\sum_{\substack{n, d \text{ odd}\\ d\sim D, n \sim N}} \frac{\chi_d(n)\rho(n)\rho'(d)}{n^{1/2+iu+w} d^{1/2-i(u+t)-w+s}}\Biggr|,
  \end{split}
\end{displaymath}
where the $D, N$-sums run over $O(\log Y_2)$ powers of 2. Recall that $\delta_0$ in $Q^{(P, Y_1)}_{1, \pm}(t, u; s, w)$ was defined in \eqref{fepara}. We combine Lemmata 3-5 to see that under the assumption \eqref{assump} we have, uniformly in $|t| \leq Y_1$ and $|u| \leq  Y_2$, 
\begin{equation}\label{3terms}
\begin{split}
 & Z(1/2 + it, 1/2+iu; \psi, \psi') \ll Y_2^{\varepsilon} \left(\frac{(TSU)^{1/4}}{\min(S, U)} + \frac{U^{3/8}}{S^A} + U^{-A}\right) +  Y_2^{\varepsilon} \sum_{P =2^{\mu}\leq Y_2^{1+\varepsilon}} \sum_{\rho, \rho'} \\
 & \quad\quad\quad  \int_{\gamma} \int_{\gamma} \min\left(\sum_{\pm}Q^{(P, Y_1)}_{1, \pm}(t, u; s, w), Q_2^{(P, Y_1, Y_2)}(t, u;  s, w) + Q_3^{(P, Y_1, Y_2)}(t, u; s, w)\right) dw\, ds \\
&\ll  Y_2^{\varepsilon} \left((Y_1Y_2)^{1/4} S^{-3/4} + Y_2^{3/8}{S^{-A}} + 1\right) + Y_2^{\varepsilon} \sum_{P =2^{\mu}\leq Y_1} \sum_{\rho, \rho', \pm}  \int_{\gamma} \int_{\gamma}  Q^{(P, Y_1)}_{1, \pm}(t, u; s, w) dw\, ds \\
&\quad\quad\quad  +  Y_2^{\varepsilon} \sum_{Y_1 \leq P =2^{\mu}\leq Y_2^{1+\varepsilon}} \sum_{\rho, \rho'}  \int_{\gamma} \int_{\gamma}  Q_2^{(P, Y_1, Y_2)}(t, u;  s, w) + Q_3^{(P, Y_1, Y_2)}(t, u; s, w) dw\, ds \\
&=: \mathcal{Q}_0 + \mathcal{Q}_1 + \mathcal{Q}_2 + \mathcal{Q}_3,
  \end{split}
\end{equation} 
say.  Here it is important to note that we may enlarge the summation ranges in the $d, n$-sums in Lemma 4 and 5 slightly to make them  independent of $t$ and $u$.   We substitute \eqref{3terms} into \eqref{toprove}. By Cauchy-Schwarz and \eqref{diri1} we have 
\begin{displaymath}
\int\int W\left(\frac{t}{Y_1}\right) W\left(\frac{u}{Y_2}\right) |\mathcal{Q}_1|^2 du\, dt  \ll \max_{\substack{P \leq Y_1,  D \leq P^{1+\varepsilon}\\ N \leq (Y_1P)^{1/2+\varepsilon}}}  Y_1Y_2^{1+\varepsilon}\left(1+\frac{N}{Y_1}\right)\left(1+\frac{D}{Y_2}\right) \ll Y_1Y_2^{1+\varepsilon}. 
\end{displaymath}
Similarly, using \eqref{diri2} instead of \eqref{diri1}, we find 
\begin{displaymath}
\begin{split}
 \int\int W\left(\frac{t}{Y_1}\right) W\left(\frac{u}{Y_2}\right) |\mathcal{Q}_2|^2 du\, dt & \ll  Y_2^{\varepsilon} \ \max_{\substack{Y_1 \leq P \leq Y_2^{1+\varepsilon} \\ D \leq  (Y_1P)^{1/2+\varepsilon}\\ N \leq (\frac{Y_1}{P})^{1/2}Y_2^{1+\varepsilon} }}  \frac{NDY_1Y_2 +ND^2 Y_1 + N^2DY_1 + (ND)^2}{DN}\\
& \ll Y_1Y_2^{1+\varepsilon}
\end{split}
\end{displaymath}
and 
\begin{displaymath}
\begin{split}
 \int\int W\left(\frac{t}{Y_1}\right) W\left(\frac{u}{Y_2}\right) |\mathcal{Q}_3|^2 du\, dt & \ll  Y_2^{\varepsilon} \ \max_{\substack{Y_1 \leq P \leq Y_2^{1+\varepsilon} \\ D \leq (\frac{Y_1}{P})^{1/2}Y_2^{1+\varepsilon}\\ N \leq  (Y_1P)^{1/2+\varepsilon} }}  \frac{NDY_1Y_2 +ND^2Y_1 + N^2DY_1 + (ND)^2}{DN}\\
 &  \ll Y_1Y_2^{1+\varepsilon}. 
\end{split}
\end{displaymath}
Finally 
we estimate trivially
\begin{displaymath}
  \int\int W\left(\frac{t}{Y_1}\right) W\left(\frac{u}{Y_2}\right) |\mathcal{Q}_0|^2 du\, dt \ll Y_2^{\varepsilon}\left(Y_1^{3/2} Y_2^{1/2} + Y_2^{3/4} Y_1 + Y_1Y_2\right)   \ll Y_1Y_2^{1+\varepsilon}
\end{displaymath}
using \eqref{assump} in both steps. The preceding four estimates establish \eqref{toprove} and complete the proof of Theorem 2.


\begin{thebibliography}{DGH}

\bibitem[Bl]{Bl} V. Blomer, \emph{On the central value of symmetric square $L$-functions}, Math. Z. \textbf{260} (2008), 755-777


  

\bibitem[DGH]{DGH} A. Diaconu, D. Goldfeld, J. Hoffstein, \emph{Multiple Dirichlet series and moments of zeta and $L$-functions},  Compositio Math.  \textbf{139}  (2003),   297--360.

\bibitem[GH]{GH} D. Goldfeld, J. Hoffstein, \emph{Eisenstein series of ${1\over 2}$-integral weight and the mean value of real Dirichlet $L$-series},  Invent. Math.  \textbf{80}  (1985),    185--208.
 
\bibitem[Go]{Go} A. Good, \emph{The square mean of Dirichlet series associated with cusp forms}, Mathematika \textbf{29} (1982),   278--295 
 

\bibitem[Ha]{Ha} G. Harcos, {\em Uniform approximate functional equation for principal $L$-functions}, Int. Math. Res. Not. \textbf{2002}, 923--932; {\em Erratum}, ibid. \textbf{2004}, 659--660

\bibitem[HM]{HM} G. Harcos, P. Michel, \emph{The subconvexity problem for Rankin-Selberg $L$-functions and equidistribution of Heegner points. II},  Invent. Math.  \textbf{163}  (2006), 581--655

\bibitem[HB]{HB} D. R. Heath-Brown, \emph{A mean value estimate for real character sums}, Acta Arith. \textbf{72} (1995),  235--275

\bibitem[HK]{HK} J. Hoffstein, A. Kontorovich, \emph{The first non-vanishing quadratic twist of an automorphic $L$-function}, preprint

\bibitem[IK]{IK} H. Iwaniec, E. Kowalski, \emph{Analytic number theory}, American Mathematical Society Colloquium Publications 53, American Mathematical Society, Providence, RI, 2004



\bibitem[Li]{Li} X. Li, \emph{Bounds for $GL(3) \times GL(2)$ $L$-functions 
and $GL(3)$ $L$-functions}, preprint


\bibitem[MV]{MV} P. Michel, A. Venkatesh, \emph{The subconvexity problem for $GL_2$}, preprint


\bibitem[Ve]{Ve} A. Venkatesh, \emph{Sparse equidistribution problems, periods bounds, and subconvexity}, Annals of Math, to appear

\bibitem[Y]{Y} M. Young, \emph{The second moment of $GL(3) \times GL(2)$ $L$-functions at special points}, preprint
  
\end{thebibliography}
\end{document}